\documentclass[12pt]{amsart}
\usepackage{epsfig}
\usepackage{young}

\theoremstyle{plain}
  \newtheorem{theorem}{Theorem}
  \newtheorem{proposition}[theorem]{Proposition}
  \newtheorem{lemma}[theorem]{Lemma}
  \newtheorem{corollary}[theorem]{Corollary}
  \newtheorem{conjecture}[theorem]{Conjecture}
\theoremstyle{definition}
  \newtheorem{definition}[theorem]{Definition}
  \newtheorem{example}[theorem]{Example}
\theoremstyle{remark}
  \newtheorem{remark}[theorem]{Remark}

\def\A{\mathcal{A}}
\def\O{\mathcal{O}}
\def\id{\mathrm{id}}
\def\flat{\mathrm{flat}}
\def\des{\mathrm{des}}
\def\x{$\times$}

\title[Bruhat order and hyperplane arrangements]
      {Bruhat order, smooth Schubert varieties, and hyperplane arrangements}
      
\author{Suho Oh, Alexander Postnikov, Hwanchul Yoo}
\address{Department of Mathematics, Massachusetts Institute of Technology, 
        77 Massachusetts Ave, Cambridge, MA 02139}
\date{September 12, 2007}
\thanks{S.O.\ was supported in part by Samsung Scholarship. 
A.P.\ was supported in part by NSF CAREER Award DMS-0504629.}

\begin{document}

\begin{abstract}
The aim of this article is to link Schubert varieties in 
the flag manifold with hyperplane arrangements.
For a permutation, we construct a certain graphical hyperplane arrangement.
We show that the generating function for regions of this arrangement
coincides with the Poincar\'e polynomial of the corresponding 
Schubert variety if and only if the Schubert variety is smooth.
We give an explicit combinatorial formula for the Poincar\'e polynomial.
Our main technical tools are chordal graphs and perfect elimination
orderings.
\end{abstract}

\maketitle

\section{Introduction}

For a permutation $w\in S_n$, let 
$P_w(q) := \sum_{u\leq w} q^{\ell(u)}$, where the sum is over
all permutations $u\in S_n$ below $w$ in the strong Bruhat order.
Geometrically, the polynomial $P_w(q)$ is the Poincar\'e polynomial
of the Schubert variety $X_w = BwB/B$ 
in the flag manifold $SL(n,\mathbb{C})/B$.

Define the {\it inversion hyperplane arrangement\/} $\A_w$ as the collection
of the hyperplanes $x_i-x_j=0$ in $\mathbb{R}^n$, for all inversions
$1\leq i < j\leq n$, $w(i)>w(j)$.  Let $R_w(q):=\sum_r q^{d(r_0,r)}$ 
be the generating function that counts regions $r$ of the arrangement $\A_w$
according to the distance $d(r_0,r)$ from the fixed initial region $r_0$.

The main result of the paper is the claim that $P_w(q)=R_w(q)$ if and 
only if the Schubert variety $X_w$ is smooth.  

According to well-known Lakshmibai-Sandhya's criterion
\cite{LS},
the Schubert variety $X_w$ is smooth if and only if the permutation $w$
avoids two patterns $3412$ and $4213$.  (Let us say that the permutation 
$w$ is smooth in this case.)  Also Carrell-Peterson \cite{C}
proved that $X_w$ is smooth if and only if the Poincar\'e polynomial
$P_w(q)$ is palindromic, that is $P_w(q) = q^{\ell(w)} \, P_w(q^{-1})$.
If $w$ is not smooth then the polynomial $P_w(q)$ is not palindromic, 
but the polynomial $R_w(q)$ is always palindromic.  So $P_w(q)\ne R_w(q)$
in this case.  On the other hand, we show that, for smooth $w$, 
the polynomials $R_w(q)$ and $P_w(q)$ satisfy the same recurrence relation. 
For the Poincar\'e polynomials $P_w(q)$, this recurrence relation was
given by Gasharov \cite{G}.  This implies that $P_w(q)=R_w(q)$ in this case.

For smooth $w$,  we present an explicit factorization 
of the polynomials $P_w(q)=R_w(q)$ as a product of $q$-numbers 
$[e_1+1]_q\cdots [e_n+1]_q$,
where $e_1,\dots,e_n$ can be computed using the left-to-right
maxima (aka records) of the permutation $w$.
In this case, the {\it inversion graph\/} $G_w$, whose edges correspond
to inversions in $w$, is a chordal graph.  The numbers $e_1,\dots,e_n$
are the roots of the chromatic polynomial $\chi_{G_w}(t)$ of the 
inversion graph.  The polynomial $\chi_{G_w}(t)$ is also the characteristic
polynomial of the inversion hyperplane arrangement $\A_w$.
We call the numbers $e_1,\dots,e_n$ the {\it exponents.}

\medskip
We thank Vic Reiner  and Jonas Sj\"ostrand 
for helpful conversations.

\section{Bruhat order and Poincar\'e polynomials}

The (strong) {\it Bruhat order\/} ``$\leq $'' on the
{\it symmetric group\/} $S_n$ is the partial order generated 
by the relations $w< w\cdot t_{ij}$ if $\ell(w)<\ell(w\cdot t_{ij})$.   
Here $t_{ij}\in S_n$ is the
transposition of $i$ and $j$; and $\ell(w)$ denotes the {\it length\/} of a
permutation $w\in S_n$, i.e., the number of inversions in $w$.

Intervals in the Bruhat order play a role in Schubert calculus and
in Kazhdan-Lusztig theory.  In this paper we concentrate on Bruhat intervals 
of the form $[\id, w]:=\{u\in S_n \mid u\leq w\}$
(where $\id\in S_n$ is the identity permutation), that is, on lower order ideals
of the Bruhat order.  They are related to {\it Schubert varieties\/}
$X_w = BwB/B$ in the flag manifold $SL(n,\mathbb{C})/B$. 
Here $B$ denotes the Borel subgroup of $SL(n,\mathbb{C})$.
The {\it Poincar\'e polynomial\/} of the Schubert variety $X_w$
is the rank generating function for the interval $[\id, w]$,
e.g., see \cite{BL}:
$$
P_w(q) = \sum_{u\leq w} q^{\ell(u)}.
$$

The well-known smoothness criterion for Schubert varieties, 
due to Lakshmibai and Sandhya,
is based on pattern avoidance. 
A permutation $w\in S_n$ contains a {\it pattern\/} $\sigma\in S_k$
if there is a subword with $k$ letters in $w$ with the same relative
order of the letters as in the permutation $\sigma$.
A permutation $w$ {\it avoids the pattern $\sigma$}
if $w$ does not contain this pattern.

\begin{theorem} {\rm (Lakshmibai-Sandhya \cite{LS})} 
For a permutation $w \in S_n$, the Schubert variety $X_w$
is smooth if and only if $w$ avoids the two patterns $3412$ and $4231$.  
\end{theorem}

We will say that $w\in S_n$ is a {\it smooth permutation\/}
if it avoids these two patterns $3412$ and $4231$.  

Another smoothness criterion, due to Carrell and Peterson, 
is given in terms of the 
Poincar\'e polynomial $P_w(q)$.  Let us say that a polynomial
$f(q)=a_0 + a_1\, q + \cdots + a_d\, q^d$ 
is {\it palindromic\/} if $f(q) = q^d f(q^{-1})$, i.e.,
$a_i = a_{d-i}$ for $i=0,\dots,d$.

\begin{theorem}
\label{th:Peterson_Criterion}
{\rm (Carrell-Peterson~\cite{C}, see also \cite[Sect.~6.2]{BL})} 
For a permutation $w \in S_n$, the Schubert variety $X_w$
is smooth if and only if the Poincar\'e polynomial $P_w(q)$
is palindromic.
\end{theorem}

\section{Inversion hyperplane arrangements}

For a graph $G$ on the vertex set $\{1,\dots,n\}$, the {\it graphical
arrangement\/} $\A_G$ is the hyperplane arrangement in $\mathbb{R}^n$ 
with hyperplanes $x_i - x_j = 0$ for all edges $(i,j)$ in $G$.
The {\it characteristic polynomial\/} $\chi_G(t)$ of the 
graphical arrangement $\A_G$ is also the {\it chromatic polynomial\/} 
of the graph $G$.  The value of $\chi_G(t)$ at a positive integer $t$
equals the number of ways to color the vertices of the graph $G$ in $t$
colors so that all neighboring pairs of vertices have different colors.
The value $(-1)^n\chi_G(-1)$ is the number of regions of $\A_G$.
The regions of $\A_G$ are in bijection with acyclic orientations 
of the graph $G$.  Recall that an {\it acyclic orientation\/} is 
a way to direct edges of $G$ so that no directed cycles are formed.
The region of $\A_G$ associated with an acyclic orientation $\O$
is described by the inequalities $x_i < x_j$ for all directed
edges $i\to j$ in $\O$.

We will study a special class of graphical arrangements.
For a permutation $w\in S_n$, the {\it inversion arrangement\/} $\A_w$ 
is the arrangement with hyperplanes $x_i - x_j = 0$ 
for each inversion $1\leq i< j\leq n$, $w(i)>w(j)$.  
Define the {\it inversion graph\/} $G_w$ as the graph on the vertex
set $\{1,\dots,n\}$ with the set of edges $\{(i,j)\mid i <j,\ 
w(i)>w(j)\}$.  The arrangement $\A_w$
is the graphical arrangement $\A_G$ for the {\it inversion graph\/} $G=G_w$.
Let $R_w$ be the number of regions in the inversion arrangement $\A_w$.

Let $B_w:=\#[\id,w]=P_w(1)$ be the number of elements in the Bruhat interval
$[\id,w]$.   Interestingly, the numbers $R_w$ and $B_w$ are 
related to each other.  

\begin{theorem}
\label{conj:R=B}
{\rm (Hultman-Linusson-Shareshian-Sj\"ostrand \cite{HLSS})} 

{\rm (1)} 
For any permutation $w\in S_n$, we have $R_w\leq B_w$.

{\rm (2)} The equality $R_w = B_w$ holds if and only if 
$w$ avoids the following four patterns
$4231$, $35142$, $42513$, $351624$.
\end{theorem}

This result was conjectured in \cite{P} and verified on a computer 
for all permutations of sizes $n\leq 8$.
This conjecture was announced as an open problem in a workshop 
in Oberwolfach in January 2007. 
A.~Hultman, S.~Linusson, J.~Shareshian, and J.~Sj\"ostrand reported 
that they proved the conjecture.  Their proof will appear on arXiv.

\begin{remark}
It was proved in \cite{P} that $R_w=B_w$ for all Grassmannian 
permutations $w$, which agrees with the 
above result.
In this case, $B_w$ counts the number of totally nonnegative
cells in the corresponding Schubert variety in the Grassmannian,
see \cite{P}.
\end{remark}

\begin{remark}  The four patterns from 
Theorem~\ref{conj:R=B} came up earlier in the literature
in at least two places.
Firstly, Gasharov and Reiner \cite{GR} showed that the Schubert variety
$X_w$ can be described by simple inclusion conditions exactly when
$w$ avoids these four patterns.
Secondly, Sj\"ostrand \cite{S} showed that the Bruhat 
interval $[\id,w]$ can be described as the set of permutations
associated with rook placements that fit inside a skew Ferrers
board if and only if $w$ avoids the same four patterns.
\end{remark}

\begin{remark}  Note that each of the four patterns from
Theorem~\ref{conj:R=B} contains one of the two patterns
from Lakshmibai-Sandhya's smoothness criterion.  Thus the theorem
implies the equality $R_w=B_w$ for all smooth permutations $w$.
\end{remark}

\section{Main results}

Let us define the $q$-analog of the number of regions
of the graphical arrangement $\A_G$,
where $G$ is a graph on the vertex set $\{1,\dots,n\}$.
For two regions $r$ and $r'$ of the arrangement $\A_G$, let $d(r,r')$ 
be the number of hyperplanes in $\A_G$ that separate $r$ and $r'$.
In other words, $d(r,r')$ is the minimal number of hyperplanes we need
to cross to go from $r$ to $r'$.
Let $r_0$ be the region of $\A_G$ that contains the point
$(1,\dots,n)$.  Define
$$
R_G(q) := \sum_{r} q^{d(r,r_0)},
$$
where the sum is over all regions $r$ of the arrangement $\A_G$.
Equivalently, the polynomial $R_G(q)$ can be described in terms
of acyclic orientations of the graph $G$.
For an acyclic orientation $\O$, let $\des(\O)$ be the number 
of edges of $G$ oriented as $i\to j$ in $\O$ where $i>j$
(descent edges).  Then 
$$
R_G(q) = \sum_\O q^{\des(\O)}, 
$$
where the sum is over all acyclic orientations $\O$ of $G$.
Indeed, for the acyclic orientation $\O$ associated with a region 
$r$ we have $\des(\O) = d(r,r_0)$.

For $w\in S_n$, let $R_w(q) := R_{G_w}(q)$ be the polynomial 
that counts the regions of 
the inversion arrangement $\A_w = A_{G_w}$. 

We are now ready to formulate the first main result of this paper.
Recall that $P_w(q) := \sum_{u\leq w} q^{\ell(u)}$ is 
the Poincar\'e polynomial of the Schubert variety.

\begin{theorem}
\label{th:P=R}
For a permutation $w\in S_n$,
we have $P_w(q) = R_w(q)$ if and only if 
$w$ is a smooth permutation, i.e., if and only if 
$w$ avoids the patterns $3412$ and $4231$. 
\end{theorem}

This result was initially conjectured during a conversation of
the second author (A.P.) with Vic Reiner.

The ``only if'' part of Theorem~\ref{th:P=R} is straightforward.
Indeed, if $w$ is not smooth then by 
Carrell-Peterson's smoothness criterion (Theorem~\ref{th:Peterson_Criterion}) 
the Poincar\'e polynomial $P_w(q)$ is not palindromic.  On the other hand,
the polynomial $R_w(q)$ is always palindromic, which follows
from the involution on the regions induced by the map
$x\mapsto -x$.  Thus $P_w(q)\ne R_w(q)$ in this case.
We will prove the ``if'' part of Theorem~\ref{th:P=R}
in Section~\ref{sec:recurrence}.

Our second result is an explicit non-recursive formula
for the polynomials $P_w(q) = R_w(q)$, when $w$ is smooth.

Let us say that an index $r\in\{1,\dots,n\}$ is a {\it record position\/} 
of a permutation $w\in S_n$ if $w(r) > \max(w(1),\dots,w(r-1))$.
The values $w(r)$ are called the {\it records\/} or 
{\it left-to-right maxima\/} of $w$.
For $i=1,\dots,n$, let $r$ and $r'$ be the record positions of $w$ 
such that $r\leq i<r'$ and there are no other record positions 
between $r$ and $r'$.  (Set $r'=+\infty$ if there are no record positions 
greater than $i$.)
Let
$$
e_i:=\#\{j\mid r\leq j<i,\ w(j)>w(i)\} +
\#\{k\mid r'\leq k\leq n,\ w(k)<w(i)\}.
$$

\begin{theorem}
\label{th:formula}
Let $w$ be a smooth permutation in $S_n$,
and let $e_1,\dots,e_n$ be the numbers constructed from $w$ as above.
Then 
$$
P_w(q) = R_w(q) = [e_1+1]_q \, [e_2+1]_q \cdots [e_n+1]_q \,.
$$
\end{theorem}

Here $[a]_q := (1-q^a)/(1-q) = 1 + q+ q^2 + \cdots + q^{a-1}$.
We will prove  Theorem~\ref{th:formula} in Section~\ref{sec:SPEO}.

\begin{example}
\label{examp:exps}
Let $w=5\,1\,6\,4\,7\,3\,2$.  
The record positions of $w$ are $1,3,5$.
We have 
$$
(e_1,\dots,e_7) = (0+3,\ 1+0,\ 0+2,\ 1+2,\ 0+0,\ 1+0,\ 2+0).
$$
Theorem~\ref{th:formula} says that 
$P_w(q)=R_w(q)  = [4]_q\, [2]_q\,[3]_q\,[4]_q\,[1]_q\,[2]_q\,[3]_q$.
\end{example}

\begin{remark}
It was known before that the Poincar\'e polynomial $P_w(q)$ for smooth $w$ 
factors as a product of $q$-numbers $[a]_q$.
Gasharov \cite{G}
(see Proposition~\ref{prop:Gasharov_factorization} below)
gave a recursive construction for such factorization.
On the other hand, Carrell gave a closed non-recursive expression 
for $P_w(q)$ as a ratio of two polynomials, see \cite{C} and
\cite[Thm.~11.1.1]{BL}.
However, it is not immediately clear from that expression that its 
denominator divides the numerator. 
One benefit of the formula in Theorem~\ref{th:formula} is that it is 
non-recursive and it involves no division.
Another combinatorial formula for $P_w(q)$ that has these features
was given by Billey, see \cite{B} and \cite[Thm.~11.1.8]{BL}. 
\end{remark}

\section{Chordal graphs and perfect elimination orderings}

A graph is called {\it chordal\/} if each of its cycles with four or more vertices has a {\it chord,} 
which is an edge joining two vertices that are not adjacent in the cycle.
A {\it perfect elimination ordering} 
in a graph $G$ is an ordering of the vertices of $G$ such that, 
for each vertex $v$ of $G$, all the neighbors of $v$ that precede $v$ in the ordering form a clique (i.e., a complete subgraph). 

\begin{theorem} {\rm (Fulkerson-Gross \cite{FG})}
A graph is chordal if and only if it has a perfect elimination ordering.
\end{theorem}

It is easy to calculate the chromatic polynomial $\chi_G(t)$ of a chordal graph $G$.
Let us pick a perfect elimination ordering $v_1,\dots,v_n$ of the vertices of $G$.  For $i=1,\dots,n$, let $e_i$ be the number of 
the neighbors of the vertex $v_i$ among the preceding vertices 
$v_{1},\dots,v_{i-1}$.  
The numbers $e_1,\dots,e_n$ are called the {\it exponents\/} of $G$.
The following formula is well-known.

\begin{proposition}
\label{prop:chromatic}
The chromatic polynomial of the chordal graph $G$ equals
$\chi_G(t) = (t-e_1)(t-e_2)\cdots (t-e_n)$.
Thus the graphical arrangement $\A_G$ has
$(-1)^n\chi_{G}(-1)=(e_1+1)(e_2+1)\cdots (e_n+1)$ regions.
\end{proposition}

For completeness sake, we include the proof, which also well-known.

\begin{proof} It is enough to prove the formula for a positive integer $t$.  Let us count the number
of coloring of vertices of $G$ in $t$ colors.  The vertex $v_1$ can be colored in $t=t-e_1$ colors.
Then the vertex $v_2$ can be colored in $t-e_2$ colors, and so on.  The vertex $v_i$ can be colored
in $t-e_i$ colors, because the $a_i$ preceding neighbors of $v_i$ already used $a_i$ {\it different\/}
colors.
\end{proof}

\begin{remark}
\label{rem:does_not_depend_on_PEO}
A chordal graph can have many different perfect elimination orderings
that lead to different sequences of exponents.  However, the
multiset (unordered sequence) $\{e_1,\dots,e_n\}$ of the exponents does not
depend on a choice of a perfect elimination order.  Indeed, by
Proposition~\ref{prop:chromatic}, the exponents $e_i$ are the roots of the
chromatic polynomial $\chi_G(t)$.  
\end{remark}

\begin{lemma}
\label{lem:R=mR}
{\rm (cf.\ Bj\"orner-Edelman-Ziegler \cite{BEZ})} \
Suppose that a graph $G$ on the vertex set $\{1,\dots,n\}$
has a vertex $v$ adjacent to $m$ vertices that satisfy the two conditions:
\begin{enumerate}
\item The set of all neighbors of $v$ is a clique in $G$.
\item 
  \begin{enumerate} 
\item All neighbors of $v$ are less than $v$, or
\item all neighbors of $v$ are greater than $v$.
  \end{enumerate}
\end{enumerate}
Then 
$
R_G(q) = [m+1]_q\, R_{G\setminus v}(q),
$
where $G\setminus v$ is the graph $G$ with the vertex $v$ removed.
\end{lemma}

This claim follows from general results of \cite{BEZ} on supersolvable
hyperplanes arrangements.  For completeness, we give a simple proof.

\begin{proof}
The polynomials $R_G(q)$ and $R_{G\setminus v}(q)$
are $\des$-generating functions for acyclic orientations
of the graphs $G$ and $G\setminus v$.

Let us fix an acyclic orientation $\O$ of the graph $G\setminus v$,
and count all ways to extend $\O$ to an acyclic orientation of $G$.  
The vertex $v$ is connected to a subset $S$ of $m$ vertices 
of the graph $G\setminus v$, which forms the clique $G|_S \simeq K_m$.
Clearly, there are $m+1$ ways to extend an acyclic orientation of 
the complete graph $K_m$ to an acyclic orientation of $K_{m+1}$.
Moreover, for each $j=0,\dots,m$, there is a unique extension of $\O$
to an acyclic orientation $\O'$ of $G$ 
such that there are exactly $j$ edges
oriented towards the vertex $v$ in $\O'$
(and $m-j$ edges oriented away from $v$).

All vertices in $S$ are less than $v$ or
all of them are greater than $v$.
In both cases we have $\sum_{\O'} q^{\des(\O')} 
= [m+1]_q\, q^{\des(\O)}$, where the sum is over extensions $\O'$ of $\O$.
Thus $R_G(q)  = [m+1]_q\, R_{G\setminus v}(q)$.
\end{proof}

\begin{definition}
\label{def:nicePEO}
For a chordal graph $G$ on the vertex set $\{1,\dots,n\}$,
we say that a perfect elimination ordering $v_1,\dots,v_n$ of the vertices 
of $G$ is {\it nice\/} if it satisfies
the following additional property.  For $i=1,\dots,n$, all neighbors of 
the vertex $v_i$ among the vertices $v_1,\dots,v_{i-1}$
are greater than $v_i$ (in the usual order on $\mathbb{Z}$), 
or all neighbors of $v_i$ among $v_1,\dots,v_{i-1}$
are less than $v_i$.
\end{definition}

For a nice perfect elimination ordering $v_1,\dots,v_n$ 
of $G$, the last vertex $v=v_n$ satisfies the conditions of 
Lemma~\ref{lem:R=mR}.  
Moreover, $v_1,\dots,v_{n-1}$ is 
a nice perfect elimination ordering of the graph 
$G\setminus v_n$.  In this case, we can inductively use 
Lemma~\ref{lem:R=mR}  to completely factor the polynomial $R_G(q)$
as $R_G(q) = [m+1]_q \,[m'+1]_q\cdots$.  
The numbers $m,m',\dots$ are exactly the exponents $e_n,e_{n-1},\dots$
(written backwards) coming from this perfect elimination ordering.

\begin{corollary}
\label{cor:nice_R}
Suppose that $G$ has a nice perfect elimination ordering of vertices.  
Let $e_1,\dots,e_n$ be the exponents of $G$.  Then we have
$$
R_G(q) = [e_1+1]_q\, [e_2+1]_q \cdots [e_n+1]_q\,.
$$
\end{corollary}


%


\section{Recurrence for polynomials $R_w(q)$}
\label{sec:recurrence}

\begin{figure}[ht]
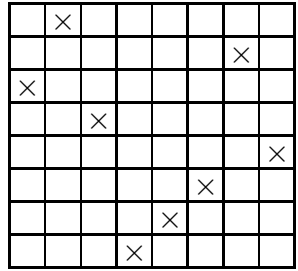

$$
\begin{Young}
&\x&&&&&&\cr
&&&&&&\x&\cr
\x&&&&&&&\cr
&&\x&&&&&\cr
&&&&&&&\x\cr
&&&&&\x&&\cr
&&&&\x&&&\cr
&&&\x&&&&\cr
\end{Young}
$$
\caption{The rook diagram $D_w$ of the permutation
$w=3\,1\,4\,8\,7\,6\,2\,5$.}
\label{fig:rook_diagram}
\end{figure}

It is convenient to represent a permutation $w\in S_n$ as
the {\it rook diagram\/} $D_w$, which the placement of $n$ non-attacking 
rooks into the boxes $(w(1),1),(w(2),2),\dots,(w(n),n)$
of the $n\times n$ board. See an example on Figure~\ref{fig:rook_diagram}.
We assume that boxes of the board
are labelled by pairs $(i,j)$ in the same way as matrix elements.
The rooks are marked by $\times$'s.

The inversion graph $G_w$ contains an edge $(i,j)$, with $i<j$,
whenever the rook in the $i$-th column of $D_w$
is located to the South-West of the rook
in the $j$-th column.  In this case, we say
that this pair of rooks {\it forms an inversion.}

Here are the rook diagrams of the two forbidden patterns 
$3412$ and $4231$ for smooth permutations:
$$
\begin{Young}
&&\x&\cr 
&&&\x\cr 
\x&&&\cr 
&\x&&\cr 
\end{Young}
\qquad
\qquad
\qquad
\begin{Young}
&&&\x\cr 
&\x&&\cr 
&&\x&\cr 
\x&&&\cr 
\end{Young}
$$
A permutation $w$ is smooth if and only if its diagram $D_w$
does not contain four rooks located in the same relative order
as in one of these diagrams $D_{3412}$ or $D_{4231}$.

\begin{figure}[ht]
\includegraphics[height=1.5in]{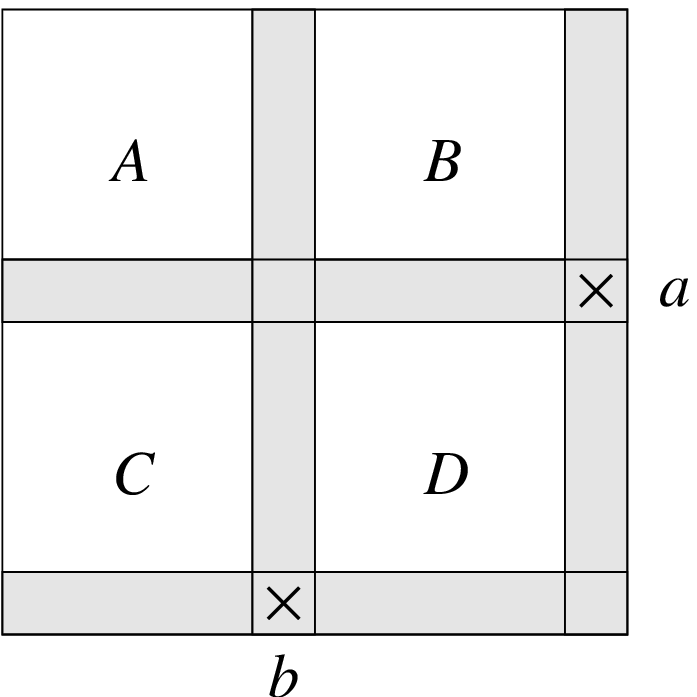}
\caption{}
\label{fig:rooks_ABCD}
\end{figure}

Let $a$ be the rook located
in the last column of $D_w$, and let $b$ be the rook located in 
the last row of $D_w$.
The row containing $a$ and the column containing $b$ subdivide
the diagram $D_w$ into the four sectors $A,B,C,D$, as shown
on Figure~\ref{fig:rooks_ABCD}.
In the case when $w(n)=n$, we assume that $a=b$ and the sectors 
$B,C,D$ are empty.

\begin{lemma}  Let $w$ be a smooth permutations.
Then its rook diagram $D_w$ has the following two properties.
{\rm (1)} Each pair of rooks located in the sector $D$
forms an inversion.
{\rm (2)} At least one of the sectors $B$ or $C$ contains
no rooks.
\end{lemma}

For example, for the rook diagram 
$D_{31487625}$ shown on Figure~\ref{fig:rook_diagram},
the sector $B$ contains one rook,
the sector $C$ contains no rooks,
and the sector $D$ contains two rooks that form an inversion.

\begin{proof}
(1) If the sector $D$ contains a pair of rooks that do not form
an inversion, then these two rooks together with the rooks $a$ and $b$
form a forbidden pattern as in the diagram $D_{4231}$.
(2) If the sector $B$ contains at least one rook and
the sector $C$ contains at least one rook, then these two rooks
together with the rooks $a$ and $b$ form a forbidden pattern
as in the diagram $D_{3412}$.
\end{proof}

Let $v_a=n$ and $v_b$ be the vertices of the inversion graph $G_w$
corresponding to the rooks $a$ and $b$.  Also let $v_1,\dots,v_k$
be the vertices of $G_w$ corresponding to the rooks inside the sector $D$.

If the sector $B$ of the rook diagram $D_w$ is empty,
then the vertex $v_b$ is connected only with the vertices
$v_1,\dots,v_k, v_a$, that form a clique in the graph $G_w$, and all these
vertices are greater
than $v_b$.  On the other hand, if the sector $C$ of the rook diagram $D_w$
is empty, then the vertex $v_a$ is connected only with the vertices
$v_b,v_1,\dots,v_k$, that form a clique, and all these vertices 
are less than $v_a$.

In both cases, the inversion graph $G_w$ satisfies the conditions
of Lemma~\ref{lem:R=mR}, where $v=v_b$ if $B$ is empty,
and $v=v_a$ if $C$ is empty.  (If both $B$ and $C$ are empty then
we can pick $v=v_a$ or $v=v_b$.)

For $w\in S_n$ and $k\in\{1,\dots,n\}$, let $w'=\flat(w,k)\in S_{n-1}$ be 
the {\it flattening\/} of the sequence $w(1),\dots,w({k-1}),
w({k+1}),\dots,w(n)$, that is, the permutation $w'$ has the 
same relative order of elements as in this sequence.
Equivalently, the rook diagram $D_{w'}$ 
is obtained from the rook diagram $D_w$ by removing 
its $k$-th column and the $w(k)$-th row.

Lemma~\ref{lem:R=mR}, together with the 
above discussion, implies the following 
recurrence relations for the polynomials $R_{w}(q)$.

\begin{proposition}
\label{prop:factor}
Let $w\in S_n$ be a smooth permutation,
and assume that $w(d)=n$ and $w(n)=e$. 
Then (at least) one of the following two
statements is true:
\begin{enumerate}
\item $w(d)>w({d+1})>\cdots > w(n)$, or
\item $w^{-1}(e)>w^{-1}({e+1})>\cdots > w^{-1}(n)$.
\end{enumerate}
In both cases, the polynomial $R_w(q)$ factors as
$$
R_w(q) = [m+1]_q \,R_{w'}(q),
$$
where $w'=\flat(w,d)$ and $m=n-d$ in case {\rm(1)},
or $w'=\flat(w,n)$ and $m=n-e$ in case {\rm(2)}.
\end{proposition}

In this proposition, case (1) means
that the sector $B$ of the rook diagram $D_w$ is empty,
and case (2) mean that the sector $C$ is empty.

Clearly, if $w$ is smooth, then the flattening $w'=\flat(w,k)$ 
is smooth as well.  The inversion graph $G_{w'}$ 
is isomorphic to the graph $G\setminus k$.
This means that, for smooth $w\in S_n$, one can inductively 
use Proposition~\ref{prop:factor} to completely factor the polynomial 
$R_w(q)$ as in Corollary~\ref{cor:nice_R}.

\begin{corollary}
\label{cor:Rwq_factor}
For a smooth permutation $w\in S_n$, 
the inversion graph $G_w$ is chordal and, moreover, it 
has a nice perfect elimination ordering.  We have
$R_w(q) = [e_1+1]_q\, [e_2+1]_q \cdots [e_n+1]_q$\,,
where $e_1,\dots,e_n$ are the exponents of the inversion
graph $G_w$.
\end{corollary}

Interestingly, Gasharov \cite{G} found exactly the same recurrence
relations for the Poincar\'e polynomials $P_w(q)$.

\begin{proposition}
\label{prop:Gasharov_factorization}
{\rm (Gasharov \cite{G}, cf.\ Lascoux \cite{L})}
The Poincar\'e polynomials $P_w(q)$, for smooth permutations $w$,
satisfy exactly the same recurrence relation as in 
Proposition~\ref{prop:factor}.
\end{proposition}

Note that Lascoux \cite{L} gave a factorization of 
the Kazhdan-Lusztig basis elements, that implies 
Proposition~\ref{prop:Gasharov_factorization}.

Propositions~\ref{prop:factor} and~\ref{prop:Gasharov_factorization},
together with the trivial claim $P_{\id}(q)=R_{\id}(q)=1$,
imply that $P_w(q)=R_w(q)$ for all smooth permutations $w$.  This finishes
the proof of Theorem~\ref{th:P=R}.

\section{Simple perfect elimination ordering}
\label{sec:SPEO}

Section~\ref{sec:recurrence} gives a recursive construction
for a nice perfect elimination ordering of the graph $G_w$, for smooth $w$. 
In this section we give a simple non-recursive construction for another
perfect elimination ordering of $G_w$.  This simple ordering 
may not be nice (see Definition~\ref{def:nicePEO}).  
However, one still can use it for calculating the 
exponents of the graph $G_w$ and factorizing the polynomials
$P_w(q)=R_w(q)$ as in Corollary~\ref{cor:Rwq_factor}.
Indeed, the multiset of the exponents does not depend
on a choice of a perfect elimination ordering
(see Remark~\ref{rem:does_not_depend_on_PEO}).

Recall that a {\it record position\/} of a permutation $w\in S_n$ is 
an index $r\in\{1,\dots,n\}$ such that 
$w(r) > \max(w(1),\dots,w(r-1))$.
Let $[a,b]$ denote the interval $\{a,a+1,\dots,b\}$ 
with the usual $\mathbb{Z}$-order of entries.

\begin{lemma}
\label{lem:speo}
For a smooth permutation $w\in S_n$
with record positions $r_1=1<r_2<\cdots<r_s$, the ordering
$$
[r_s,n],\
[r_{s-1},\,r_s-1],\dots,[r_2,\,r_3-1],\ [r_1,\,r_{2}-1]
$$
of the set $\{1,\dots,n\}$
is a perfect elimination ordering of the inversion graph $G_w$.
\end{lemma}

\begin{example} (cf.\ Example~\ref{examp:exps})
The permutation $w=5\,1\,6\,4\,7\,3\,2$
has records $5,6,7$ and record positions $1,3,5$.
Lemma~\ref{lem:speo} says that the ordering
$\underline{5,\,6,\,7},\ \underline{3,\,4},\ \underline{1,\,2}$
is a perfect elimination ordering of the inversion graph $G_w$.
Figure~\ref{fig:inv_graph} displays this inversion graph $G_w$.
For each vertex $i=1,\dots,7$ of $G_w$, we wrote $i$ inside a circle
and $w(i)$ below it.
The exponents of this graph
(i.e., the numbers of edges going to the left
from the vertices) are $0,1,2,2,3,3,1$.
\begin{figure}[ht]
\includegraphics[height=1.2in]{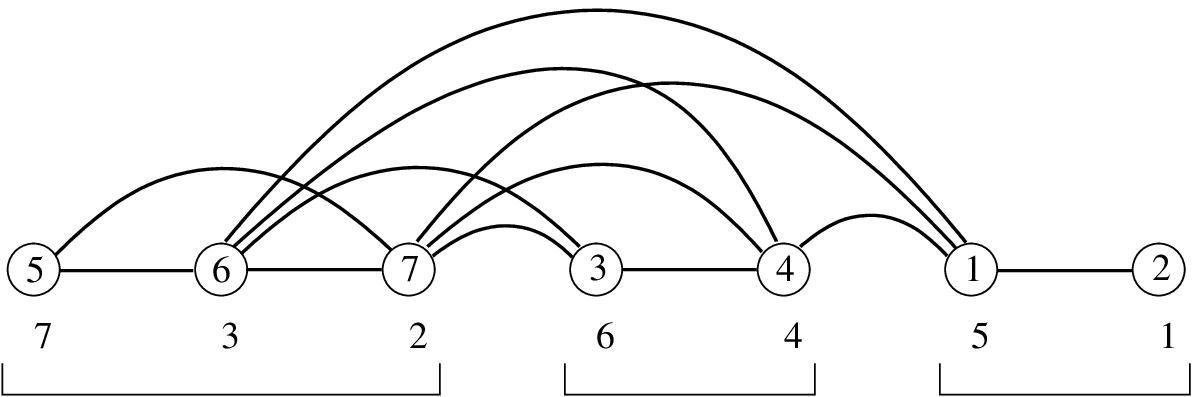}
\caption{}
\label{fig:inv_graph}
\end{figure}
\end{example}

\begin{proof}[Proof of Lemma~\ref{lem:speo}]
Suppose that this ordering of vertices of $G_w$ is not a perfect elimination 
ordering.  This means that there is a vertex $i$ connected in $G_w$ with 
vertices $j$ and $k$, preceding $i$ in the order, such that 
the vertices $j$ and $k$ are not connected by an edge in $G_w$.
Let us consider three cases.

I.  The vertices $i,j,k$ belong to the same interval 
$I_p:=[r_p,r_{p+1}-1]$, for some $p\in\{1,\dots,s\}$.  
(Here we assume that $r_{s+1} = n+1$.) We have $k<j<i$ and 
$w(k)>w(i)$, $w(j)>w(i)$, but  $w(k)<w(j)$,
because $(k,i)$ and $(j,k)$ are edges of $G_w$ but $(k,j)$ is not an edge.
The value $w(r_p)$ is the maximal value of $w$ on the interval $I_p$.  
Since $w(k)<w(j)$ is not the maximal value of $w$ on $I_p$,
we have $r_p\ne k$ and so $r_p<k$.  
Thus $r_p<k<j<i$ and the values $w(r_p), w(k), w(j), w(i)$ form a forbidden
$4231$ pattern in $w$.  So $w$ is not smooth.  Contradiction.

II.  The vertices $i,j$ are in the same interval $I_p$
and the vertex $k$ belongs to a different interval $I_q$.
Then $q>p$, because the vertex $k$ precedes $i$ in the order.
In this case we have $j<i<k$, $w(j)>w(i)$, $w(i)>w(k)$.
This implies that $w(j)>w(k)$ that is $(j,k)$ is an edge
in the inversion graph $G_w$.  Contradiction.

III.  The vertex $i$ belongs to the interval $I_p$ and 
the vertices $j,k$ do not belong to $I_p$.  Assume that $j<k$
and that $j$ belongs to $I_q$.  Then $q>p$.
In this case, $i<j<k$, $w(i)>w(j)$, $w(i)>w(k)$, and $w(j)<w(k)$.
The record value $w(r_q)$ is greater than $w(i)$.
This implies that $w(r_q)>w(i)>w(j)$.  In particular, $w(r_q)\ne w(j)$ and,
thus, $r_q\ne j$.  We have $i<r_q<j<k$ and the values
$w(i), w(r_q), w(j), w(k)$ form a forbidden $3412$ pattern.  Contradiction.
\end{proof}

\begin{proof}[Proof of Theorem~\ref{th:formula}]
Let us calculate the exponents of the inversion graph $G_w$
for a smooth permutation $w\in S_n$ using the perfect elimination
ordering from Lemma~\ref{lem:speo}.
Suppose that $i\in I_p$.   Then the exponent $e_i$ of the vertex $i$
equals the number of neighbors of the vertex $i$ in the graph $G_w$
among the preceding vertices, that is among the vertices in the sets
$\{r_p,\dots,i-1\}$ and $I_{p+1}\cup I_{p+2}\cup \dots$.
In other words, the exponent $e_i$ equals 
$$
\#\{j\mid r_p\leq j< i,\ w(j)> w(i)\} + 
\#\{k\mid k\geq r_{p+1},\ w(k)<w(i)\}.
$$
This is exactly the expression for $e_i$ from Theorem~\ref{th:formula}.
The result follows from Corollary~\ref{cor:Rwq_factor}.
\end{proof}

\section{Final remarks}

Our proof of Theorem~\ref{th:P=R} is based on a recurrence relation.
It would be interesting to give more direct combinatorial proof of 
Theorem~\ref{th:P=R} based on a bijection between elements of the Bruhat 
interval $[\id,w]$ and regions of the arrangement $\A_w$.

\medskip

It would be interesting to better understand the relationship
between Bruhat intervals $[\id,w]$ and the hyperplane arrangement $\A_w$.
One can construct a directed graph $\Gamma_w$ on the regions of $\A_w$.
Two regions $r$ and $r'$ are connected by a directed edge $(r,r')$
if these two regions are adjacent (i.e., separated by a single hyperplane)
and $r$ is more close to $r_0$ than $r'$.
For example, for the longest permutation $w_0$, the graph $\Gamma_{w_0}$
is the Hasse diagram of the {\it weak\/} Bruhat order.
It is true that, for any smooth permutation $w\in S_n$,
the graph $\Gamma_w$ is isomorphic to a subgraph of the
Hasse diagram of the Bruhat interval $[\id,w]$?

\medskip

It would be interesting to explain Theorem~\ref{th:P=R}
from a geometrical point of view.  Is it possible to link 
the arrangement $\A_w$ and the polynomial $R_w(q)$ with 
the cohomology ring of of the Schubert variety $X_w$?
Is it possible to define a related ring 
structure on the regions of $\A_w$?

\medskip

The statement of Theorem~\ref{th:P=R} can extended to any
finite Weyl group $W$, as follows.
For a Weyl group element $w\in W$, let
$P_w(q) := \sum_u q^{\ell(w)}$, where the sum is over
all $u\in W$ such that $u\leq w$ in the Bruhat order on $W$.
Define the arrangement $\A_w$ as the collection of hyperplanes
$\alpha(x)=0$ for all roots $\alpha$ in the corresponding root
system such that $\alpha>0$ and $w(\alpha)<0$.  Let $r_0$ 
be the region of $\A_w$ that contains the fundamental chamber
of the corresponding Coxeter arrangement.
Define $R_w(q):= \sum_r q^{d(r_0,r)}$, where the sum
is over all regions of the arrangement $\A_w$
and $d(r_0,r)$ is the number of hyperplanes separating
$r_0$ and $r$.
Let $X_w = BwB/B$ be the Schubert variety in
the corresponding generalized flag manifold $G/B$.
Details about (rational) smoothness of Schubert varieties
$X_w$ can be found in \cite{BL}.

\begin{conjecture}  The equality $P_w(q)=R_w(q)$ holds
if and only if the Schubert variety $X_w$ is rationally smooth.
\end{conjecture}

Finally, let us mention that the inverse of  
Corollary~\ref{cor:nice_R} might be true.

\begin{conjecture}
For a graph $G$,
the polynomial $R_G(q)$ can be factorized as a product
of $q$-numbers if and only if the graph $G$ has 
a nice perfect elimination order.
\end{conjecture}


\begin{thebibliography}{HLSS}

\bibitem[B]{B} S.~Billey: Pattern avoidance and rational smoothness of
Schubert varieties,
{\it Adv.\ in Math.\ \bf 139} (1998), 141--156.

\bibitem[BL]{BL} 
S.~Billey, V.~Lakshmibai: 
{\it Singular Loci of Schubert Varieties,}
Progress in Mathematics, Vol.~182, Birkh\"auser, Boston, 2000.

\bibitem[BEZ]{BEZ} A.~Bj\"orner, P.~H.~Edelman, G.~M.~Ziegler:
Hyperplane arrangements with a lattice of regions,
{\it Discrete Comput.\ Geom.\ \bf 5} (1990), no.~3, 263--288.


\bibitem[C]{C} 
J.~B.~Carrell: The Bruhat graph of a Coxeter group, 
a conjecture of Deodhar, and rational smoothness of 
Schubert varieties, 
{\it Proceedings of Symposia in Pure Math.\ \bf 56} (1994), 53--61.

\bibitem[FG]{FG} D.~R.~Fulkerson, O.~A.~Gross: Incidence matrices and interval graphs, {\it Pacific J.\ Math \bf 15}, 
835--855.

\bibitem[G]{G} V.~Gasharov: Factoring the Poincare polynomials for the Bruhat
order on $S_n$, {\it Journal of Combinatorial Theory, Series A \bf 83} (1998), 
159--164.

\bibitem[GR]{GR} V.~Gasharov, V.~Reiner:
Cohomology of smooth Schubert varieties in partial flag manifolds. 
{\it J.\ London Math.\ Soc.} (2) {\bf 66} (2002), no.~3, 550--562. 

\bibitem[HLSS]{HLSS} A.~Hultman, S.~Linusson, J.~Shareshian, J.~Sj\"ostrand:
From Bruhat intervals to intersections lattices and a conjecture of Postnikov,
in preparation.

\bibitem[L]{L} 
A.~Lascoux: 
Ordonner la groupe sym\'etrique: pourquoi utiliser
l'alg\`ebre de Iwahori-Hecke,
{\it Proc.\ ICM Berlin,} Doc.\ Math.\ 1998 Extra Vol.\ III, 
355--364 (electronic).

\bibitem[LS]{LS} V.~Lakshmibai, B.~Sandhya: Criterion for smoothness of Schubert
varieties in $SL(n)/B$, {\it Proceedings of the Indian Academy of Science
(Mathematical Sciences) \bf 100} (1990), 45--52. MR 91c:14061.

\bibitem[P]{P} A.~Postnikov: Total positivity, Grassmannians, and networks, 
arXiv: math/ 0609764v1 [math.CO].

\bibitem[S]{S} 
J.~Sj\"ostrand:
Bruhat intervals as rooks on skew Ferrers boards, 
{\it J.\ Combin.\ Theory Ser.\ A} (7) 114 (2007), 1182--1198.





\end{thebibliography}
\end{document}